\documentclass{article}
\usepackage{amssymb}
\usepackage{amsfonts}
\usepackage{amsmath}

\newtheorem{theorem}{Theorem}

\newtheorem{corollary}[theorem]{Corollary}

\newtheorem{definition}[theorem]{Definition}

\newtheorem{proposition}[theorem]{Proposition}
\newtheorem{remark}[theorem]{Remark}

\newenvironment{proof}[1][Proof]{\noindent\textbf{#1.} }{\ \rule{0.5em}{0.5em}}
\input{tcilatex}
\begin{document}

\title{Polynomials Related to Harmonic Numbers and Evaluation of Harmonic
Number Series II}
\author{Ayhan Dil and Veli Kurt \\
Department of Mathematics, Akdeniz University, 07058 Antalya Turkey\\
adil@akdeniz.edu.tr, vkurt@akdeniz.edu.tr}
\maketitle

\begin{abstract}
In this paper we focus on $r-$geometric polynomials, $r-$exponential
polynomials and their harmonic versions. It is shown that harmonic versions
of these polynomials and their generalizations are useful to obtain closed
forms of some series related to harmonic numbers.


\end{abstract}

\section{Introduction}

\qquad In \cite{DK1} the concept of $harmonic-$ geometric polynomials and $%
harmonic-$ exponential polynomials are introduced and $hyperharmonic$
generalizations of these polynomials and numbers are obtained. Furthermore
it is shown that these polynomials are quite useful to obtain closed forms
of some series related to harmonic numbers. In this paper, we extend this
analysis to $r-versions$ of these polynomials and numbers.

Boyadzhiev\ $\cite{B}$ has presented\ and discussed the following
transformation formula:%
\begin{equation}
\sum_{n=0}^{\infty }\frac{g^{\left( n\right) }\left( 0\right) }{n!}f\left(
n\right) x^{n}=\sum_{n=0}^{\infty }\frac{f^{\left( n\right) }\left( 0\right) 
}{n!}\sum_{k=0}^{n}\QATOPD\{ \} {n}{k}x^{k}g^{\left( k\right) }\left(
x\right)   \label{1}
\end{equation}%
where $f$, $g$ are appropriate functions and $\QATOPD\{ \} {n}{k}$ are
Stirling numbers of the second kind.

One of the principal objectives of the present paper is to give closed forms
of some series related to harmonic numbers as well. To this end, we give a
useful generalization of $\left( \ref{1}\right) $ which contains $r-$%
Stirling numbers of the second kind instead of Stirling numbers of the
second kind as:%
\begin{equation}
\sum_{n=r}^{\infty }\frac{g^{\left( n\right) }\left( 0\right) }{n!}\binom{n}{%
r}\frac{r!}{n^{r}}f_{r}\left( n\right) x^{n}=\sum_{n=r}^{\infty }\frac{%
f^{\left( n\right) }\left( 0\right) }{n!}\sum_{k=0}^{n}\QATOPD\{ \}
{n}{k}_{r}x^{k}g^{\left( k\right) }\left( x\right) \text{,}  \label{mf}
\end{equation}%
where $f_{r}\left( x\right) $ denotes the Maclaurin series of $f\left(
x\right) $ exclude the first $r-1$ terms.

Thanks to formula $\left( \ref{mf}\right) $ we introduce the concept of $r-$
geometric and $r-$ exponential polynomials and numbers. We obtain explicit
relations between the $r-$versions and the classical versions of these
polynomials and numbers. Besides, we present harmonic (and hyperharmonic)
versions of $r-$ geometric and $r-$ exponential polynomials and numbers as
well. The short lists of all these polynomials and numbers are given.

On the other hand formula $\left( \ref{mf}\right) $ and harmonic $r-$
geometric polynomials enables us to obtain closed forms of the following
series%
\begin{equation}
\sum_{n=r}^{\infty }\binom{n}{r}r!n^{m-r}H_{n}x^{n}\text{,}  \label{mfr}
\end{equation}%
where $m$ and $r$ are integers such that $m\geq r$ and $H_{n}$ is the $n$-th
partial sum of the harmonic series

In the rest of this section we introduce some important notions.

\textbf{Stirling numbers of the first and second kind}

Stirling numbers of the first kind $\QATOPD[ ] {n}{k}$ and Stirling numbers
of the second kind $\QATOPD\{ \} {n}{k}$ are quite important in
combinatorics $\left( \cite{BG, BQ, C, Ri}\right) $. Briefly for the
integers $n\geq k\geq 0;$ $\QATOPD[ ] {n}{k}$ counts the number of
permutations of $n$ elements with exactly $k$ cycles and $\QATOPD\{ \} {n}{k}
$ counts the number of ways to partition a set with $n$ elements into $k$
disjoint, nonempty subsets $\left( \cite{C}\right) $.

We note that for $n\geq k\geq 1$, the following identity holds for Stirling
numbers of the second kind%
\begin{equation}
\QATOPD\{ \} {n}{k}=\QATOPD\{ \} {n-1}{k-1}+k\QATOPD\{ \} {n-1}{k}.
\label{4}
\end{equation}

There is a certain generalization of these numbers namely $r$-Stirling
numbers $\left( \cite{Br}\right) $ which are similar to the weighted
Stirling numbers $\left( \cite{CA1, CA2}\right) $. Represantation and
combinatorial meanings of these numbers are as follows $\left( \cite{Br}%
\right) $:

$r-$Stirling numbers of the first kind;%
\begin{eqnarray*}
\QATOPD[ ] {n}{k}_{r} &=&\text{The number of permutations of the set }%
\left\{ 1,2,\ldots ,n\right\} \text{ having} \\
&&k\text{ cycles, such that the numbers }1,2,...,r\text{ are in distinct
cylecs,}
\end{eqnarray*}

$r-$Stirling numbers of the second kind;%
\begin{eqnarray*}
\QATOPD\{ \} {n}{k}_{r} &=&\text{ The number of partitions of the set }%
\left\{ 1,2,\ldots ,n\right\} \text{ into } \\
&&k\text{ non-empty disjoint subsets, such that the numbers } \\
&&1,2,...,r\text{ are in distinct subsets.}
\end{eqnarray*}

Specializing $r=0$ gives the classical Stirling numbers.

The $r-$Stirling numbers of the second kind satisfy the same recurrence
relation as $\left( \ref{4}\right) $ except for the initial conditions, i.e $%
\left( \cite{Br}\right) $%
\begin{eqnarray}
\QATOPD\{ \} {n}{k}_{r} &=&0\text{, \ \ }n<r,  \notag \\
\QATOPD\{ \} {n}{k}_{r} &=&\delta _{k,r}\text{, \ \ }n=r,  \label{5} \\
\QATOPD\{ \} {n}{k}_{r} &=&\QATOPD\{ \} {n-1}{k-1}_{r}+k\QATOPD\{ \}
{n-1}{k}_{r}\text{, \ \ }n>r.  \notag
\end{eqnarray}%
\textbf{Exponential polynomials and numbers}

Exponential polynomials (or single variable Bell polynomials) $\phi
_{n}\left( x\right) $ are defined by$\left( \cite{BL1, B2, G, Ri}\right) $%
\begin{equation}
\phi _{n}\left( x\right) :=\sum_{k=0}^{n}\QATOPD\{ \} {n}{k}x^{k}.  \label{6}
\end{equation}

The first few exponential polynomials are:%
\begin{equation}
\begin{tabular}{|l|}
\hline
$\phi _{0}\left( x\right) =1\text{,}$ \\ \hline
$\phi _{1}\left( x\right) =x\text{,}$ \\ \hline
$\phi _{2}\left( x\right) =x+x^{2}\text{,}$ \\ \hline
$\phi _{3}\left( x\right) =x+3x^{2}+x^{3}\text{,}$ \\ \hline
$\phi _{4}\left( x\right) =x+7x^{2}+6x^{3}+x^{4}\text{.}$ \\ \hline
\end{tabular}
\label{7}
\end{equation}

The well known exponential numbers (or Bell numbers)\ are obtained by
setting $x=1$ in $\phi _{n}\left( x\right) $, i.e $\left( \cite{BL2, C, CG}%
\right) $%
\begin{equation}
\phi _{n}:=\phi _{n}\left( 1\right) =\sum_{k=0}^{n}\QATOPD\{ \} {n}{k}.
\label{8}
\end{equation}%
Hence the first few exponential numbers are:%
\begin{equation}
\phi _{0}=1\text{, }\phi _{1}=1\text{, }\phi _{2}=2\text{, }\phi _{3}=5\text{%
, }\phi _{4}=15\text{.}  \label{9}
\end{equation}

In \cite{DK} the authors obtained some fundemental properties of the
exponential polynomials and numbers using Euler-Seidel matrices method as:%
\begin{equation}
\phi _{n+1}\left( x\right) =x\sum_{k=0}^{n}\binom{n}{k}\phi _{k}\left(
x\right)  \label{10}
\end{equation}%
and%
\begin{equation}
\phi _{n+1}=\sum_{k=0}^{n}\binom{n}{k}\phi _{k}.  \label{11}
\end{equation}

Recently, Mezo $\left( \cite{M}\right) $ has defined the "$r-$Bell
polynomials and numbers" as:%
\begin{equation}
B_{n,r}\left( x\right) =\sum_{k=0}^{n}\QATOPD\{ \} {n+r}{k+r}_{r}x^{k}
\label{12}
\end{equation}%
and%
\begin{equation}
\text{ }B_{n,r}=\sum_{k=0}^{n}\QATOPD\{ \} {n+r}{k+r}_{r}  \label{12+}
\end{equation}%
respectively. $r-$Exponential polynomials and numbers which we discuss in
the present paper are slightly different than the $r-$Bell polynomials and
numbers in \cite{M}.

\textbf{Geometric polynomials and numbers}

Geometric polynomials are defined in \cite{B, S, ST} as follows:%
\begin{equation}
F_{n}\left( x\right) :=\sum_{k=0}^{n}\QATOPD\{ \} {n}{k}k!x^{k}.  \label{13}
\end{equation}

The first few geometric polynomials are:%
\begin{equation}
\begin{tabular}{|l|}
\hline
$F_{0}\left( x\right) =1\text{,}$ \\ \hline
$F_{1}\left( x\right) =x\text{,}$ \\ \hline
$F_{2}\left( x\right) =x+2x^{2}\text{,}$ \\ \hline
$F_{3}\left( x\right) =x+6x^{2}+6x^{3}\text{,}$ \\ \hline
$F_{4}\left( x\right) =x+14x^{2}+36x^{3}+24x^{4}\text{.}$ \\ \hline
\end{tabular}
\label{14}
\end{equation}

Specializing $x=1$ in $\left( \ref{13}\right) $ we get geometric numbers (or
ordered Bell numbers) $F_{n}$ as $\left( \cite{B, ST, W}\right) $:%
\begin{equation}
F_{n}:=F_{n}\left( 1\right) =\sum_{k=0}^{n}\QATOPD\{ \} {n}{k}k!.
\label{i14}
\end{equation}

The first few geometric numbers are:%
\begin{equation}
F_{0}=1\text{, }F_{1}=1\text{, }F_{2}=3\text{, }F_{3}=13\text{, }F_{4}=75%
\text{.}  \label{i15}
\end{equation}

Boyadzhiev $\left( \cite{B}\right) $\ introduced the "general geometric
polynomials" as%
\begin{equation}
F_{n,r}\left( x\right) =\frac{1}{\Gamma \left( r\right) }\sum_{k=0}^{n}%
\QATOPD\{ \} {n}{k}\Gamma \left( k+r\right) x^{k}\text{,}  \label{16+}
\end{equation}%
where Re$\left( r\right) >0$. In the fifth section we will deal with the
general geometric polynomials.

Exponential and geometric polynomials are connected by the following
integral relation $\left( \cite{B}\right) $%
\begin{equation}
F_{n}\left( z\right) =\int_{0}^{\infty }\phi _{n}\left( z\lambda \right)
e^{-\lambda }d\lambda .  \label{16}
\end{equation}

In \cite{DK} the authors also obtained some fundemental properties of the
geometric polynomials and numbers using Euler- Seidel matrices method as:%
\begin{equation}
F_{n+1}\left( x\right) =x\sum_{k=0}^{n}\binom{n+1}{k}F_{k}\left( x\right)
\label{15}
\end{equation}%
and%
\begin{equation}
F_{n}=\sum_{k=0}^{n-1}\binom{n}{k}F_{k}.  \label{15+}
\end{equation}

By means of $r-$Stirling numbers Nyul $\left( \cite{GN}\right) $ introduced "%
$r-$geometric polynomials and numbers (or $r-$Fubini or ordered $r-$Bell
polynomials and numbers)" are respectively as follows:%
\begin{equation}
F_{n,r}\left( x\right) =\sum_{k=0}^{n}\left( k+r\right) !\QATOPD\{ \}
{n+r}{k+r}_{r}x^{k}  \label{17}
\end{equation}%
and%
\begin{equation*}
F_{n,r}=\sum_{k=0}^{n}\left( k+r\right) !\QATOPD\{ \} {n+r}{k+r}_{r}.
\end{equation*}

In this work $r-$geometric polynomials come up naturally in an application
of the generalized transformation formula as well.

Notice that, our concept of $r-$geometric polynomials is slightly different
than in \cite{GN}.

\textbf{Harmonic and Hyperharmonic numbers}

The $n$-th harmonic number is the $n$-th partial sum of the harmonic series: 
\begin{equation}
H_{n}:=\sum_{k=1}^{n}\frac{1}{k}\text{,}  \label{18}
\end{equation}

where $H_{0}=0.$

For an integer $\alpha >1$, let

\begin{equation}
H_{n}^{(\alpha )}:=\sum_{k=1}^{n}H_{k}^{(\alpha -1)}\text{,}  \label{i18}
\end{equation}%
with $H_{n}^{(1)}:=H_{n}$, be the $n$-th hyperharmonic number of order $%
\alpha $ $\left( \cite{BG, CG}\right) $.

These numbers can be expressed in terms of binomial coefficients and
ordinary harmonic numbers as $\left( \cite{CG, MD}\right) $:%
\begin{equation}
H_{n}^{(\alpha )}=\binom{n+\alpha -1}{\alpha -1}(H_{n+\alpha -1}-H_{\alpha
-1}).  \label{20}
\end{equation}

Well-known generating functions of the harmonic and hyperharmonic numbers
are given by%
\begin{equation}
\sum_{n=1}^{\infty }H_{n}x^{n}=-\frac{\ln \left( 1-x\right) }{1-x}
\label{21}
\end{equation}%
and%
\begin{equation}
\sum_{n=1}^{\infty }H_{n}^{\left( \alpha \right) }x^{n}=-\frac{\ln \left(
1-x\right) }{\left( 1-x\right) ^{\alpha }}.  \label{22}
\end{equation}%
respectively $\left( \cite{DM}\right) $.

The following relations connect harmonic and hyperharmonic numbers with the
Stirling and $r-$Stirling numbers of the first kind $\left( \cite{BG}\right) 
$:%
\begin{equation}
\QATOPD[ ] {k+1}{2}=k!H_{k}\text{,}  \label{23}
\end{equation}%
and%
\begin{equation}
k!H_{k}^{\left( r\right) }=\QATOPD[ ] {n+r}{r+1}_{r}\text{.}  \label{24}
\end{equation}

\section{Generalization of the transformation formula}

\qquad In this section firstly we mention Boyadzhiev's Theorem $4.1$ in \cite%
{B} and give a useful generalization of it. As a result of this
generalization we introduce $r-$geometric polynomials and numbers.

Suppose we are given an entire function $f$\ and a function $g$, analytic in
a region containing the annulus $K=\{z:r<|z|<R\}$ where $0<r<R$. Hence these
functions have following series expansions,%
\begin{equation*}
f\left( x\right) =\sum_{n=0}^{\infty }p_{n}x^{n}\text{ and }g\left( x\right)
=\sum_{n=-\infty }^{\infty }q_{n}x^{n}.
\end{equation*}%
Now we are ready to state Boyadzhiev's theorem.

\begin{theorem}
\label{Main Theorem}$\left( \cite{B}\right) $ Let the functions $f$ and $g$
be described as above. If the series%
\begin{equation*}
\sum_{n=-\infty }^{\infty }q_{n}f\left( n\right) x^{n}
\end{equation*}%
converges absolutely on $K$, then%
\begin{equation}
\sum_{n=-\infty }^{\infty }q_{n}f\left( n\right) x^{n}=\sum_{m=0}^{\infty
}p_{m}\sum_{k=0}^{m}\QATOPD\{ \} {m}{k}x^{k}g^{\left( k\right) }\left(
x\right)  \label{0}
\end{equation}%
holds for all $x\in K$.
\end{theorem}

Stirling numbers of the second kind appear in the formula $\left( \ref{0}%
\right) $ due to $\left( xD\right) $ operator. Our aim is to get a more
general formula than $\left( \ref{0}\right) $\ which contains $r-$Stirling
numbers of the second kind instead of Stirling numbers of the second kind.
Accordingly, first we generalize the operator $\left( xD\right) $.

\subsection{Generalization of the\textbf{\ operator} $\left( xD\right) $}

The operator $\left( xD\right) $ operates a function $f\left( x\right) $ as;%
\begin{equation}
\left( xD\right) f\left( x\right) :=xf^{\prime }\left( x\right) \text{,}
\label{o1}
\end{equation}%
where $f^{\prime }$ is the first derivative of the function $f.$

For any $m$-times differentiable function $f$ we have $\left( \cite{B}%
\right) $,%
\begin{equation}
\left( xD\right) ^{m}f\left( x\right) =\sum_{k=0}^{m}\QATOPD\{ \}
{m}{k}x^{k}f^{\left( k\right) }\left( x\right) .  \label{o2}
\end{equation}%
This fact can be easily proven with induction on $m$ by the help of $\left( %
\ref{4}\right) $.

Our first aim is to generalize the operator $\left( xD\right) $.\ Later we
use this generalization to obtain $r-$ geometric and $r-$exponential
polynomials and numbers. In the light of this motivation and after a plenty
of observations we arrive the following definition.

\begin{definition}
\label{GO}Let $f$ be a function which is at least $m-$times differentiable
and $r$ be a nonnegative integer. Then the action $\left( xD_{r}\right) $ is%
\begin{equation}
\left( xD_{r}\right) ^{m}f\left( x\right) :=\left\{ 
\begin{array}{cc}
0 & ,m<r \\ 
\left( xD\right) ^{m-r}x^{r}f^{\left( r\right) }\left( x\right) & ,m\geq r%
\end{array}%
\right. .  \label{g1}
\end{equation}
\end{definition}

An equivalent statement of this definition is given by the following
proposition.

\begin{proposition}
By applying $\left( xD_{r}\right) $ $m-$times to a function $f$ which is at
least $m-$times differentiable, then the following \ equation holds%
\begin{equation}
\left( xD_{r}\right) ^{m}f\left( x\right) =\sum_{k=0}^{m}\QATOPD\{ \}
{m}{k}_{r}x^{k}f^{\left( k\right) }\left( x\right) \text{,}  \label{g2}
\end{equation}%
where $m\geq r$.
\end{proposition}

\begin{proof}
It follows from induction on $m$, in the light of Definition $\ref{GO}$\ and
recurrence relation $\left( \ref{5}\right) $.
\end{proof}

The equation $\left( \ref{g2}\right) $\ is a generalization of the equation $%
\left( \ref{o2}\right) $\ since setting $r=0$ in $\left( \ref{g2}\right) $\
gives the equation $\left( \ref{o2}\right) $.

\begin{corollary}
Let $n$ be an integer, then%
\begin{equation}
\left( xD_{r}\right) ^{m}x^{n}=n^{m-r}\binom{n}{r}r!x^{n}.  \label{g3}
\end{equation}
\end{corollary}

\subsection{Generalization of the transformation formula}

\qquad Now we give our main theorem that is a generalization of Theorem 1.

\begin{theorem}
\label{MT}Let $f\left( x\right) $ be an entire function and $g\left(
x\right) $ be an analytic function on the annulus $K=\{z,s<|z|<S\}$, where $%
0\leq s<S$. Suppose that their power series be given as%
\begin{equation*}
f\left( x\right) =\sum_{m=0}^{\infty }p_{m}x^{m}\text{ and }g\left( x\right)
=\sum_{n=-\infty }^{\infty }q_{n}x^{n}.
\end{equation*}%
If the series 
\begin{equation}
\sum_{n=-\infty }^{\infty }q_{n}\binom{n}{r}\frac{r!}{n^{r}}f_{r}\left(
n\right) x^{n}  \label{ac}
\end{equation}
where $r$ is a nonnegative integer and $f_{r}\left( x\right) $ denotes the
power series $\sum_{m=r}^{\infty }p_{m}x^{m}$, converges absolutely on $K$,
then%
\begin{equation}
\sum_{n=-\infty }^{\infty }q_{n}\binom{n}{r}\frac{r!}{n^{r}}f_{r}\left(
n\right) x^{n}=\sum_{m=r}^{\infty }p_{m}\sum_{k=0}^{m}\QATOPD\{ \}
{m}{k}_{r}x^{k}g^{\left( k\right) }\left( x\right)   \label{g4}
\end{equation}%
holds for all $x\in K$.
\end{theorem}

\begin{proof}
By considering the power series expansion of $g\left( x\right) $ with $%
\left( \ref{g2}\right) $ and $\left( \ref{g3}\right) $ we have%
\begin{equation}
\sum_{n=-\infty }^{\infty }q_{n}\binom{n}{r}n^{m-r}r!x^{n}=\sum_{k=0}^{m}%
\QATOPD\{ \} {m}{k}_{r}x^{k}g^{\left( k\right) }\left( x\right)   \label{g5}
\end{equation}%
where $m$ and $r$\ are integer such that $m\geq r\geq 0$. If we multiply
both sides of the equation $\left( \ref{g5}\right) $ by $p_{m}$ and sum on $m
$ from $r$ to infinity we get%
\begin{equation*}
\sum_{n=-\infty }^{\infty }q_{n}\binom{n}{r}\frac{r!}{n^{r}}%
\sum_{m=r}^{\infty }p_{m}n^{m}x^{n}=\sum_{m=r}^{\infty
}p_{m}\sum_{k=0}^{m}\QATOPD\{ \} {m}{k}_{r}x^{k}g^{\left( k\right) }\left(
x\right) \text{,}
\end{equation*}%
since $\left( \ref{ac}\right) $\ is converges absolutely on $K$. This
completes proof.
\end{proof}

\begin{corollary}
\label{gsv}Let $g$ be an analytic function on the disk $D=\{z,0\leq |z|<S\}$
then%
\begin{equation}
\sum_{n=r}^{\infty }\frac{g^{\left( n\right) }\left( 0\right) }{n!}\binom{n}{%
r}\frac{r!}{n^{r}}f_{r}\left( n\right) x^{n}=\sum_{n=r}^{\infty }\frac{%
f^{\left( n\right) }\left( 0\right) }{n!}\sum_{k=0}^{n}\QATOPD\{ \}
{n}{k}_{r}x^{k}g^{\left( k\right) }\left( x\right) .  \label{g6}
\end{equation}
\end{corollary}

Most of the results in the subsequent sections depend on the Corollary $\ref%
{gsv}$.

\begin{remark}
Specializing $r=0$ in the Theorem $\ref{MT}$\ we turn back to the Theorem $%
4.1$ of Boyadzhiev $\left( \cite{B}\right) $. Therefore from now on we are
interested in the case $r\geq 1$.
\end{remark}

\section{$r-$exponential and $r-$geometric polynomials and numbers}

\qquad Stirling numbers of the first and second kind are notable in many
branches of mathematics, especially in combinatorics, computational
mathematics and computer sciences $\left( \cite{AS, BQ, C, CG, GKP}\right) $%
. Importance of the exponential polynomials and numbers are substantially
because of their direct connection with Striling numbers. $r-$Stirling
numbers $\left( \cite{Br}\right) $ are one of the reputable generalizations
of Stirling numbers. Therefore introduction of the concepts of the $r-$
exponential and $r-$ geometric polynomials and numbers are good motivation
for us.

\subsection{$r-$exponential polynomials and numbers}

Firstly we consider $g\left( x\right) =e^{x}$ in the equation $\left( \ref%
{g6}\right) $. Hence we get%
\begin{equation}
\sum_{n=r}^{\infty }\binom{n}{r}\frac{r!}{n^{r}}f_{r}\left( n\right) \frac{%
x^{n}}{n!}=e^{x}\sum_{n=r}^{\infty }\frac{f^{\left( n\right) }\left(
0\right) }{n!}\sum_{k=0}^{n}\QATOPD\{ \} {n}{k}_{r}x^{k}.  \label{g7}
\end{equation}%
The finite sum on the RHS is a generalization of exponential polynomials. We
call these polynomials as "$r-$exponential polynomials" and indicate them
with $_{r}\phi _{n}\left( x\right) $. Hence%
\begin{equation}
_{r}\phi _{n}\left( x\right) :=\sum_{k=0}^{n}\QATOPD\{ \} {n}{k}_{r}x^{k}.
\label{g8}
\end{equation}

The first few $r-$exponential polynomials are:

\begin{equation}
\begin{tabular}{|l|l|l|l|}
\hline
$_{r}\phi _{n}\left( x\right) $ & $r=1$ & $r=2$ & $r=3$ \\ \hline
$n=0$ & $0$ & $0$ & $0$ \\ \hline
$n=1$ & $x$ & $0$ & $0$ \\ \hline
$n=2$ & $x+x^{2}$ & $x^{2}$ & $0$ \\ \hline
$n=3$ & $x+3x^{2}+x^{3}$ & $2x^{2}+x^{3}$ & $x^{3}$ \\ \hline
$n=4$ & $x+7x^{2}+6x^{3}+x^{4}$ & $4x^{2}+5x^{3}+x^{4}$ & $3x^{3}+x^{4}$ \\ 
\hline
\end{tabular}
\label{Lrep}
\end{equation}

Similar to the classical case, "$r-$exponential numbers" can be defined by
setting $x=1$ in $\left( \ref{g8}\right) $\ i.e,%
\begin{equation}
_{r}\phi _{n}:=\sum_{k=0}^{n}\QATOPD\{ \} {n}{k}_{r}.  \label{g8+}
\end{equation}

Hence the first few $r-$exponential numbers are:

\begin{equation}
\begin{tabular}{|l|l|l|l|}
\hline
$_{r}\phi _{n}$ & $r=1$ & $r=2$ & $r=3$ \\ \hline
$n=0$ & $0$ & $0$ & $0$ \\ \hline
$n=1$ & $1$ & $0$ & $0$ \\ \hline
$n=2$ & $2$ & $1$ & $0$ \\ \hline
$n=3$ & $5$ & $3$ & $1$ \\ \hline
$n=4$ & $15$ & $10$ & $4$ \\ \hline
\end{tabular}
\label{Lren}
\end{equation}

Now we give an explicit formula which connects $r-$exponential polynomials
with the classical exponential polynomials. Also this formula allows us to
calculate $_{r}\phi _{n}\left( x\right) $ easily.

\begin{proposition}
We have%
\begin{equation}
_{r}\phi _{n+r}\left( x\right) =x^{r}\sum_{k=0}^{n}\binom{n}{k}r^{n-k}\phi
_{k}\left( x\right) \text{,}  \label{g9}
\end{equation}%
where $n$ and $r$ are nonnegative integers.
\end{proposition}

\begin{proof}
Let $m$ be an integer such that $m\geq r\geq 0$ and we specialize $f\left(
x\right) =x^{m}$ in $\left( \ref{g7}\right) $.Then we get%
\begin{equation*}
_{r}\phi _{m}\left( x\right) e^{x}=\sum_{n=r}^{\infty }\binom{n}{r}\frac{r!}{%
n!}n^{m-r}x^{n}.
\end{equation*}%
RHS of this equation can be written as%
\begin{equation*}
x^{r}\sum_{n=0}^{\infty }\left( n+r\right) ^{m-r}\frac{x^{n}}{n!}%
=x^{r}\sum_{k=0}^{m-r}\binom{m-r}{k}r^{m-r-k}\sum_{n=0}^{\infty }n^{k}\frac{%
x^{n}}{n!}.
\end{equation*}%
Considering the definition of the operator $\left( xD\right) $ this becomes%
\begin{equation*}
x^{r}\sum_{k=0}^{m-r}\binom{m-r}{k}r^{m-r-k}\left( xD\right) ^{k}e^{x}.
\end{equation*}%
The equation $\left( \ref{o2}\right) $ enables us to write%
\begin{equation*}
x^{r}e^{x}\sum_{k=0}^{m-r}\binom{m-r}{k}r^{m-r-k}\phi _{k}\left( x\right) .
\end{equation*}%
Comparision of the LHS and the RHS completes the proof.
\end{proof}

Similar relation can be given between classical exponential numbers and $r-$%
exponential numbers as a corollary.

\begin{corollary}
\begin{equation}
_{r}\phi _{n+r}=\sum_{k=0}^{n}\binom{n}{k}r^{n-k}\phi _{k}.  \label{g10}
\end{equation}
\end{corollary}

The following corollary shows that the equation $\left( \ref{g9}\right) $ is
a generalization of the equation $\left( \ref{10}\right) $.

\begin{corollary}
\begin{equation*}
\phi _{n+1}\left( x\right) =x\sum_{k=0}^{n}\binom{n}{k}\phi _{k}\left(
x\right) \text{.}
\end{equation*}
\end{corollary}

\subsection{$r-$geometric polynomials and numbers}

By considering $g\left( x\right) =\frac{1}{1-x}$ in the equation $\left( \ref%
{g6}\right) $ we get%
\begin{equation}
\sum_{n=r}^{\infty }\binom{n}{r}\frac{r!}{n^{r}}f_{r}\left( n\right) x^{n}=%
\frac{1}{1-x}\sum_{n=r}^{\infty }\frac{f^{\left( n\right) }\left( 0\right) }{%
n!}\sum_{k=0}^{n}\QATOPD\{ \} {n}{k}_{r}k!\left( \frac{x}{1-x}\right) ^{k}.
\label{rFg}
\end{equation}%
We call the finite sum of the RHS as "$r-$geometric polynomials" and
indicate them with $_{r}F_{n}\left( x\right) $. Hence%
\begin{equation}
_{r}F_{n}\left( x\right) :=\sum_{k=0}^{n}\QATOPD\{ \} {n}{k}_{r}k!x^{k}.
\label{g11}
\end{equation}

The first few $r-$geometric polynomials are:

\begin{equation}
\begin{tabular}{|l|l|l|l|}
\hline
$_{r}F_{n}\left( x\right) $ & $r=1$ & $r=2$ & $r=3$ \\ \hline
$n=0$ & $0$ & $0$ & $0$ \\ \hline
$n=1$ & $x$ & $0$ & $0$ \\ \hline
$n=2$ & $x+2x^{2}$ & $2x^{2}$ & $0$ \\ \hline
$n=3$ & $x+6x^{2}+6x^{3}$ & $4x^{2}+6x^{3}$ & $6x^{3}$ \\ \hline
$n=4$ & $x+14x^{2}+36x^{3}+24x^{4}$ & $8x^{2}+30x^{3}+24x^{4}$ & $%
18x^{3}+24x^{4}$ \\ \hline
\end{tabular}
\label{Lrgp}
\end{equation}

We define "$r-$geometric numbers" by specializing $x=1$ in $\left( \ref{g11}%
\right) $\ as%
\begin{equation}
_{r}F_{n}:=\sum_{k=0}^{n}\QATOPD\{ \} {n}{k}_{r}k!.  \label{g12}
\end{equation}

The first few $r-$geometric numbers are:

\begin{equation}
\begin{tabular}{|l|l|l|l|}
\hline
$_{r}F_{n}$ & $r=1$ & $r=2$ & $r=3$ \\ \hline
$n=0$ & $0$ & $0$ & $0$ \\ \hline
$n=1$ & $1$ & $0$ & $0$ \\ \hline
$n=2$ & $3$ & $2$ & $0$ \\ \hline
$n=3$ & $13$ & $10$ & $6$ \\ \hline
$n=4$ & $75$ & $62$ & $42$ \\ \hline
\end{tabular}
\label{Lrgn}
\end{equation}

The following proposition gives an explicit formula between $r-$geometric
polynomials and generalized geometric polynomials which have given by the
equation $\left( \ref{16+}\right) $.

\begin{proposition}
For any nonnegative integers $n$ and $r$ we have%
\begin{equation}
_{r}F_{n+r}\left( x\right) =x^{r}r!\sum_{k=0}^{n}\binom{n}{k}%
r^{n-k}F_{k,r+1}\left( x\right)  \label{g13}
\end{equation}
\end{proposition}

\begin{proof}
Let $m$ be a nonnegative integer such that $m\geq r$. By setting $f\left(
x\right) =x^{m}$ in $\left( \ref{rFg}\right) $ we get%
\begin{equation*}
\frac{1}{1-x}\text{ }_{r}F_{m}\left( \frac{x}{1-x}\right)
=\sum_{n=r}^{\infty }\binom{n}{r}r!n^{m-r}x^{n}\text{.}
\end{equation*}%
Rearranging RHS gives%
\begin{equation*}
x^{r}r!\sum_{k=0}^{m-r}\binom{m-r}{k}r^{m-r-k}\sum_{n=0}^{\infty }\binom{n+r%
}{r}n^{k}x^{n}\text{.}
\end{equation*}%
We can write this by means of $\left( xD\right) $\ operator as%
\begin{equation*}
x^{r}r!\sum_{k=0}^{m-r}\binom{m-r}{k}r^{m-r-k}\left( xD\right) ^{k}\frac{1}{%
\left( 1-x\right) ^{r+1}}.
\end{equation*}%
Considering the fact that (equation $\left( 3.26\right) $ in \cite{B})%
\begin{equation*}
\left( xD\right) ^{k}\frac{1}{\left( 1-x\right) ^{r+1}}=\frac{1}{\left(
1-x\right) ^{r+1}}F_{k,r+1}\left( \frac{x}{1-x}\right) 
\end{equation*}%
completes the proof.
\end{proof}

A similar result between numbers is as follows.

\begin{corollary}
\begin{equation}
_{r}F_{n+r}=r!\sum_{k=0}^{n}\binom{n}{k}r^{n-k}F_{k,r+1}  \label{g14}
\end{equation}
\end{corollary}

Owing to $\left( \ref{g14}\right) $,\ we give the following relations for
classical geometric polynomials and numbers as a corollary.

\begin{corollary}
\begin{equation}
F_{n+1}\left( x\right) =x\sum_{k=0}^{n}\binom{n}{k}F_{k,2}\left( x\right) 
\text{,}  \label{g15}
\end{equation}%
\begin{equation}
F_{n+1}=\sum_{k=0}^{n}\binom{n}{k}F_{k,2}.  \label{g16}
\end{equation}
\end{corollary}

\section{Harmonic $r-$geometric and harmonic $r-$exponential polynomials and
numbers}

\qquad We introduce the concepts of $harmonic-$geometric and $harmonic-$%
exponential polynomials and numbers in $\cite{DK1}$. Along this section we
follow similar approach in $\left( \cite{DK1}\right) $\ to investigate
harmonic $r-$geometric and harmonic $r-$exponential polynomials and numbers.

\subsection{Harmonic $r-$geometric polynomials and numbers}

\qquad We consider the generating function of harmonic numbers\ as the
function $g$\ in the transformation formula $\left( \ref{g6}\right) .$ From 
\cite{DK1} we have

\begin{equation}
g^{\left( k\right) }\left( z\right) =\frac{k!\left( H_{k}-\ln \left(
1-z\right) \right) }{\left( 1-z\right) ^{k+1}}  \label{r1}
\end{equation}%
and

\begin{equation}
g^{\left( k\right) }\left( 0\right) =k!H_{k}.  \label{r2}
\end{equation}

With the help of Theorem $\ref{MT}$ we state the following transformation
formula for harmonic numbers.

\begin{proposition}
Let $r$ be a nonnegative integer and $f$ be an entire function. Then we have%
\begin{eqnarray}
&&\sum_{n=r}^{\infty }\binom{n}{r}H_{n}\frac{r!}{n^{r}}f_{r}\left( n\right)
x^{n}  \notag \\
&=&\frac{1}{1-x}\sum_{n=r}^{\infty }\frac{f^{\left( n\right) }\left(
0\right) }{n!}\sum_{k=0}^{n}\QATOPD\{ \} {n}{k}_{r}k!H_{k}\left( \frac{x}{1-x%
}\right) ^{k}  \label{r3} \\
&&-\frac{\ln \left( 1-x\right) }{1-x}\sum_{n=r}^{\infty }\frac{f^{\left(
n\right) }\left( 0\right) }{n!}\sum_{k=0}^{n}\QATOPD\{ \} {n}{k}_{r}k!\left( 
\frac{x}{1-x}\right) ^{k}.  \notag
\end{eqnarray}
\end{proposition}

\begin{proof}
Employing $\left( \ref{r1}\right) $ and $\left( \ref{r2}\right) $ in $\left( %
\ref{g6}\right) $ gives the statement.
\end{proof}

Second part of the RHS of the equation $\left( \ref{r3}\right) $ contains $r-
$geometric polynomials which are familiar to us from the previous section.
But the first part contains a new family of polynomials which is a
generalization of $harmonic-$geometric polynomials $\left( \cite{DK1}\right) 
$. We call them as "harmonic $r-$geometric polynomials and indicate them
with $_{r}F_{n}^{h}\left( x\right) $. Thus

\begin{equation}
_{r}F_{n}^{h}\left( x\right) :=\sum_{k=0}^{n}\QATOPD\{ \}
{n}{k}_{r}k!H_{k}x^{k}.  \label{r4}
\end{equation}

The first few harmonic $r-$geometric polynomials are:

\begin{equation}
\begin{tabular}{|l|l|l|l|}
\hline
$_{r}F_{n}^{h}\left( x\right) $ & $r=1$ & $r=2$ & $r=3$ \\ \hline
$n=0$ & $0$ & $0$ & $0$ \\ \hline
$n=1$ & $x$ & $0$ & $0$ \\ \hline
$n=2$ & $x+3x^{2}$ & $3x^{2}$ & $0$ \\ \hline
$n=3$ & $x+9x^{2}+11x^{3}$ & $6x^{2}+11x^{3}$ & $11x^{3}$ \\ \hline
$n=4$ & $x+21x^{2}+66x^{3}+50x^{4}$ & $12x^{2}+55x^{3}+50x^{4}$ & $%
33x^{3}+50x^{4}$ \\ \hline
\end{tabular}
\label{Lhrgp}
\end{equation}

"Harmonic $r-$geometric numbers" can be defined by setting $x=1$ in $\left( %
\ref{r4}\right) $, i.e%
\begin{equation}
_{r}F_{n}^{h}:=\sum_{k=0}^{n}\QATOPD\{ \} {n}{k}_{r}k!H_{k}.  \label{r4+}
\end{equation}

The first few harmonic $r-$geometric numbers are:

\begin{equation}
\begin{tabular}{|l|l|l|l|}
\hline
$_{r}F_{n}^{h}$ & $r=1$ & $r=2$ & $r=3$ \\ \hline
$n=0$ & $0$ & $0$ & $0$ \\ \hline
$n=1$ & $1$ & $0$ & $0$ \\ \hline
$n=2$ & $4$ & $3$ & $0$ \\ \hline
$n=3$ & $21$ & $17$ & $11$ \\ \hline
$n=4$ & $138$ & $117$ & $83$ \\ \hline
\end{tabular}
\label{Lhrgn}
\end{equation}

Hence with this notation we state the formula $\left( \ref{r3}\right) $
simply as%
\begin{eqnarray}
&&\sum_{n=r}^{\infty }\binom{n}{r}H_{n}\frac{r!}{n^{r}}f_{r}\left( n\right)
x^{n}  \notag \\
&=&\frac{1}{1-x}\sum_{n=r}^{\infty }\frac{f^{\left( n\right) }\left(
0\right) }{n!}\left\{ _{r}F_{n}^{h}\left( \frac{x}{1-x}\right)
-_{r}F_{n}\left( \frac{x}{1-x}\right) \ln \left( 1-x\right) \right\} .
\label{r5}
\end{eqnarray}

Due to the following corollary we obtain closed forms of some series related
to harmonic numbers and binomial coefficients.

\begin{corollary}
\begin{equation}
\sum_{n=r}^{\infty }\binom{n}{r}r!n^{m-r}H_{n}x^{n}=\frac{1}{1-x}\left\{
_{r}F_{m}^{h}\left( \frac{x}{1-x}\right) -_{r}F_{m}\left( \frac{x}{1-x}%
\right) \ln \left( 1-x\right) \right\} \text{,}  \label{r6}
\end{equation}%
where $m$ and $r$ are integers such that $m\geq r$.
\end{corollary}

\begin{proof}
It follows by setting $f\left( x\right) =x^{m}$ in the equation $\left( \ref%
{r5}\right) $.
\end{proof}

\begin{remark}
Formula $\left( \ref{r6}\right) $ allow us to calculate closed forms of
several harmonic number series. The case $r=1$ in $\left( \ref{r6}\right) $\
has been analyzed in \cite{DK1} already.
\end{remark}

The case $r=2$ gives%
\begin{equation}
\sum_{n=2}^{\infty }n^{m-1}\left( n-1\right) H_{n}x^{n}=\frac{1}{1-x}\left\{
_{2}F_{m}^{h}\left( \frac{x}{1-x}\right) -_{2}F_{m}\left( \frac{x}{1-x}%
\right) \ln \left( 1-x\right) \right\} .  \label{gs1}
\end{equation}%
Hence some series and their closed forms that we get from $\left( \ref{gs1}%
\right) $\ are as follows:

For $m=2$ we have%
\begin{equation}
\sum_{n=2}^{\infty }n\left( n-1\right) H_{n}x^{n}=\frac{x^{2}\left\{ 3-2\ln
\left( 1-x\right) \right\} }{\left( 1-x\right) ^{3}}\text{.}  \label{gs1a}
\end{equation}

For $m=3$ we have%
\begin{equation}
\sum_{n=2}^{\infty }n^{2}\left( n-1\right) H_{n}x^{n}=\frac{x^{2}\left\{
6+5x-\left( 4+2x\right) \ln \left( 1-x\right) \right\} }{\left( 1-x\right)
^{4}}\text{,}  \label{gs1b}
\end{equation}%
and so on.

The case $r=3$ gives%
\begin{eqnarray}
&&\sum_{n=3}^{\infty }n^{m-2}\left( n-1\right) \left( n-2\right) H_{n}x^{n} 
\notag \\
&=&\frac{1}{1-x}\left\{ _{3}F_{m}^{h}\left( \frac{x}{1-x}\right)
-_{3}F_{m}\left( \frac{x}{1-x}\right) \ln \left( 1-x\right) \right\} \text{.}
\label{gs2}
\end{eqnarray}%
Hence some series and their closed forms that we get from $\left( \ref{gs2}%
\right) $\ are as follows:

For $m=3$ we have%
\begin{equation}
\sum_{n=3}^{\infty }n\left( n-1\right) \left( n-2\right) H_{n}x^{n}=\frac{%
x^{3}\left\{ 11-6\ln \left( 1-x\right) \right\} }{\left( 1-x\right) ^{4}}%
\text{.}  \label{gs3}
\end{equation}

For $m=4$ we have%
\begin{equation}
\sum_{n=3}^{\infty }n^{2}\left( n-1\right) \left( n-2\right) H_{n}x^{n}=%
\frac{x^{3}\left\{ 33+17x-\left( 18+6x\right) \ln \left( 1-x\right) \right\} 
}{\left( 1-x\right) ^{5}}\text{,}  \label{gs4}
\end{equation}%
and so on.

Now we give a summation formula for the multiple series.

\begin{proposition}
\begin{eqnarray}
&&\sum_{n=r}^{\infty }\left( \sum_{k=0}^{n-r}\binom{k}{r}\binom{n+s-k}{s}%
r!k^{m-r}H_{k}\right) x^{n}  \notag \\
&=&\sum_{n=r}^{\infty }\left( \sum_{0\leq k_{1}\leq k_{2}\leq \cdots \leq
k_{s+1}\leq n}\binom{k_{1}}{r}r!k_{1}^{m-r}H_{k_{1}}\right) x^{n}
\label{gs5} \\
&=&\frac{1}{\left( 1-x\right) ^{s+2}}\left\{ _{r}F_{m}^{h}\left( \frac{x}{1-x%
}\right) -_{r}F_{m}\left( \frac{x}{1-x}\right) \ln \left( 1-x\right) \right\}
\notag
\end{eqnarray}
\end{proposition}

\begin{proof}
Multiplying both sides of the equation $\left( \ref{r6}\right) $ with the
Newton binomial series and considering that%
\begin{equation*}
\sum_{k=0}^{n-r}\binom{k}{r}\binom{n+s-k}{s}r!k^{m-r}H_{k}=\sum_{0\leq
k_{1}\leq k_{2}\leq \cdots \leq k_{s+1}\leq n}\binom{k_{1}}{r}%
r!k_{1}^{m-r}H_{k_{1}}
\end{equation*}%
we get the statement.
\end{proof}

By setting $r=2$ and $s=0$ in the formula $\left( \ref{gs5}\right) $ we can
give the following applications:

For $m=2$ we have%
\begin{equation}
\sum_{n=2}^{\infty }\left( \sum_{k=2}^{n}k\left( k-1\right) H_{k}\right)
x^{n}=\frac{x^{2}\left\{ 3-2\ln \left( 1-x\right) \right\} }{\left(
1-x\right) ^{4}}\text{.}  \label{gs5a}
\end{equation}

For $m=3$ we have%
\begin{equation}
\sum_{n=2}^{\infty }\left( \sum_{k=2}^{n}k^{2}\left( k-1\right) H_{k}\right)
x^{n}=\frac{x^{2}\left\{ 6+5x-\left( 4+2x\right) \ln \left( 1-x\right)
\right\} }{\left( 1-x\right) ^{5}}\text{,}  \label{gs5b}
\end{equation}%
and so on.

\begin{remark}
By the help of $\left( \ref{23}\right) $ we can state $_{r}F_{n}^{h}\left(
x\right) $ and $_{r}F_{n}^{h}$ in terms of $r-$Stirling numbers of the
second kind and Stirling numbers of the first kind as%
\begin{equation}
_{r}F_{n}^{h}\left( x\right) =\sum_{k=0}^{n}\QATOPD\{ \} {n}{k}_{r}\QATOPD[ ]
{k+1}{2}x^{k}  \label{r7}
\end{equation}%
and%
\begin{equation}
_{r}F_{n}^{h}=\sum_{k=0}^{n}\QATOPD\{ \} {n}{k}_{r}\QATOPD[ ] {k+1}{2}.
\label{r8}
\end{equation}
\end{remark}

\subsection{Harmonic $r-$exponential polynomials and numbers}

\qquad Bearing in mind the similarity of exponential and geometric
polynomials and being inspried by the definition of harmonic exponential
polynomials and numbers we arrive the following definition.

\begin{definition}
\label{hrb}For the nonnegative integers $n$ and $r,$ "harmonic $r-$%
exponential polynomials" and "harmonic $r-$exponential numbers" are defined
respectively as%
\begin{equation}
_{r}\phi _{n}^{h}\left( x\right) :=\sum_{k=0}^{n}\QATOPD\{ \}
{n}{k}_{r}H_{k}x^{k}  \label{r9}
\end{equation}%
and%
\begin{equation}
_{r}\phi _{n}^{h}:=\sum_{k=0}^{n}\QATOPD\{ \} {n}{k}_{r}H_{k}.  \label{r10}
\end{equation}
\end{definition}

The first few harmonic $r-$exponential polynomials are:

\begin{equation}
\begin{tabular}{|l|l|l|l|}
\hline
$_{r}\phi _{n}^{h}\left( x\right) $ & $r=1$ & $r=2$ & $r=3$ \\ \hline
$n=0$ & $0$ & $0$ & $0$ \\ \hline
$n=1$ & $x$ & $0$ & $0$ \\ \hline
$n=2$ & $x+\frac{3}{2}x^{2}$ & $\frac{3}{2}x^{2}$ & $0$ \\ \hline
$n=3$ & $x+\frac{9}{2}x^{2}+\frac{11}{6}x^{3}$ & $3x^{2}+\frac{11}{6}x^{3}$
& $\frac{11}{6}x^{3}$ \\ \hline
$n=4$ & $x+\frac{21}{2}x^{2}+11x^{3}+\frac{25}{12}x^{4}$ & $6x^{2}+\frac{55}{%
6}x^{3}+\frac{25}{12}x^{4}$ & $\frac{11}{2}x^{3}+\frac{25}{12}x^{4}$ \\ 
\hline
\end{tabular}
\label{Lhrep}
\end{equation}

The first few harmonic $r-$exponential numbers are:

\begin{equation}
\begin{tabular}{|l|l|l|l|}
\hline
$_{r}\phi _{n}^{h}$ & $r=1$ & $r=2$ & $r=3$ \\ \hline
$n=0$ & $0$ & $0$ & $0$ \\ \hline
$n=1$ & $1$ & $0$ & $0$ \\ \hline
$n=2$ & $\frac{5}{2}$ & $\frac{3}{2}$ & $0$ \\ \hline
$n=3$ & $\frac{22}{3}$ & $\frac{29}{6}$ & $\frac{11}{6}$ \\ \hline
$n=4$ & $\frac{295}{12}$ & $\frac{69}{4}$ & $\frac{91}{12}$ \\ \hline
\end{tabular}
\label{Lhren}
\end{equation}

\begin{remark}
Definition $\ref{hrb}$ enables us to extend the relation $\left( \ref{16}%
\right) \ $as%
\begin{equation}
_{r}F_{n}^{h}\left( z\right) =\int_{0}^{\infty }\text{ }_{r}\phi
_{n}^{h}\left( z\lambda \right) e^{-\lambda }d\lambda .  \label{r11}
\end{equation}
\end{remark}

\section{Hyperharmonic $r-$geometric and hyperharmonic $r-$exponential
polynomials and numbers}

\qquad For the completeness of this work, now we consider hyperharmonic
numbers and their transformations. In this way we could generalize almost
all results of \cite{DK1} and in previous sections of the present paper.

\subsection{Hyperharmonic $r-$geometric polynomials and numbers}

\qquad Similar to the previous section, let us consider the function $g$ in
the transformation formula $\left( \ref{g6}\right) $\ as the generating
function of the hyperharmonic numbers. From \cite{DK1} we have

\begin{equation}
g^{\left( k\right) }\left( x\right) =\frac{\Gamma \left( k+\alpha \right) }{%
\Gamma \left( \alpha \right) }\frac{1}{\left( 1-z\right) ^{\alpha +k}}\left(
H_{k+\alpha -1}-H_{\alpha -1}-\ln \left( 1-x\right) \right)  \label{r12}
\end{equation}%
and

\begin{equation}
g^{\left( k\right) }\left( 0\right) =k!H_{k}^{\left( \alpha \right) }.
\label{r13}
\end{equation}%
Now we give a transformation formula for hyperharmonic numbers.

\begin{proposition}
\label{phht}For integers $r\geq 0$ and $\alpha \geq 1$ we have%
\begin{eqnarray}
&&\sum_{n=r}^{\infty }\binom{n}{r}H_{n}^{\left( \alpha \right) }\frac{r!}{%
n^{r}}f_{r}\left( n\right) x^{n}  \notag \\
&=&\frac{1}{\left( 1-z\right) ^{\alpha }}\sum_{n=r}^{\infty }\frac{f^{\left(
n\right) }\left( 0\right) }{n!}\sum_{k=0}^{n}\QATOPD\{ \}
{n}{k}_{r}k!H_{k}^{\left( \alpha \right) }\left( \frac{x}{1-x}\right) ^{k}
\label{r14} \\
&&-\frac{\ln \left( 1-x\right) }{\left( 1-z\right) ^{\alpha }}%
\sum_{n=r}^{\infty }\frac{f^{\left( n\right) }\left( 0\right) }{n!}\frac{1}{%
\Gamma \left( \alpha \right) }\sum_{k=0}^{n}\QATOPD\{ \} {n}{k}_{r}\Gamma
\left( k+\alpha \right) \left( \frac{x}{1-x}\right) ^{k}.  \notag
\end{eqnarray}
\end{proposition}

\begin{proof}
Consideration $\left( \ref{r12}\right) $ and $\left( \ref{r13}\right) $ in $%
\left( \ref{g6}\right) $ give the statement.
\end{proof}

The first part of the RHS is a generalization of harmonic $r-$geometric
polynomials which contains hyperharmonic numbers instead of harmonic
numbers. We call these polynomials as "hyperharmonic $r-$geometric
polynomials" and indicate them with $_{r}F_{n,\alpha }^{h}\left( x\right) $.
Thus%
\begin{equation}
_{r}F_{n,\alpha }^{h}\left( x\right) =\sum_{k=0}^{n}\QATOPD\{ \}
{n}{k}_{r}k!H_{k}^{\left( \alpha \right) }x^{k}  \label{r15}
\end{equation}

The first few hyperharmonic $r-$geometric polynomials are:

Case $\alpha =2$

\begin{equation}
\begin{tabular}{|l|l|l|}
\hline
$_{r}F_{n,2}^{h}\left( x\right) $ & $r=1$ & $r=2$ \\ \hline
$n=0$ & $0$ & $0$ \\ \hline
$n=1$ & $x$ & $0$ \\ \hline
$n=2$ & $x+5x^{2}$ & $5x^{2}$ \\ \hline
$n=3$ & $x+15x^{2}+26x^{3}$ & $10x^{2}+26x^{3}$ \\ \hline
$n=4$ & $x+35x^{2}+156x^{3}+154x^{4}$ & $20x^{2}+130x^{3}+154x^{4}$ \\ \hline
\end{tabular}
\label{hhrgp2}
\end{equation}

Case $\alpha =3$

\begin{equation}
\begin{tabular}{|l|l|l|}
\hline
$_{r}F_{n,3}^{h}\left( x\right) $ & $r=1$ & $r=2$ \\ \hline
$n=0$ & $0$ & $0$ \\ \hline
$n=1$ & $x$ & $0$ \\ \hline
$n=2$ & $x+7x^{2}$ & $7x^{2}$ \\ \hline
$n=3$ & $x+21x^{2}+47x^{3}$ & $14x^{2}+47x^{3}$ \\ \hline
$n=4$ & $x+49x^{2}+282x^{3}+342x^{4}$ & $28x^{2}+235x^{3}+342x^{4}$ \\ \hline
\end{tabular}
\label{hhrgp3}
\end{equation}

The second part of the RHS of $\left( \ref{r14}\right) $\ contains also a
generalization of the polynomials, "general geometric polynomials",\ which
we mention with the equation $\left( \ref{16+}\right) $. We call these
polynomials as "general $r-$geometric polynomials" and indicate them with $%
_{r}F_{n,\alpha }\left( x\right) $. Hence%
\begin{equation}
_{r}F_{n,\alpha }\left( x\right) =\frac{1}{\Gamma \left( \alpha \right) }%
\sum_{k=0}^{n}\QATOPD\{ \} {n}{k}_{r}\Gamma \left( k+\alpha \right) x^{k}.
\label{r16}
\end{equation}

The first few general $r-$geometric polynomials are:

Case $\alpha =2$

\begin{equation}
\begin{tabular}{|l|l|l|}
\hline
$_{r}F_{n,2}\left( x\right) $ & $r=1$ & $r=2$ \\ \hline
$n=0$ & $0$ & $0$ \\ \hline
$n=1$ & $2x$ & $0$ \\ \hline
$n=2$ & $2x+6x^{2}$ & $6x^{2}$ \\ \hline
$n=3$ & $2x+18x^{2}+24x^{3}$ & $12x^{2}+24x^{3}$ \\ \hline
$n=4$ & $2x+42x^{2}+144x^{3}+120x^{4}$ & $24x^{2}+120x^{3}+120x^{4}$ \\ 
\hline
\end{tabular}
\label{grgp2}
\end{equation}%
and

Case $\alpha =3$

\begin{equation}
\begin{tabular}{|l|l|l|}
\hline
$_{r}F_{n,3}\left( x\right) $ & $r=1$ & $r=2$ \\ \hline
$n=0$ & $0$ & $0$ \\ \hline
$n=1$ & $3x$ & $0$ \\ \hline
$n=2$ & $3x+12x^{2}$ & $12x^{2}$ \\ \hline
$n=3$ & $3x+36x^{2}+60x^{3}$ & $24x^{2}+60x^{3}$ \\ \hline
$n=4$ & $3x+84x^{2}+360x^{3}+360x^{4}$ & $48x^{2}+300x^{3}+360x^{4}$ \\ 
\hline
\end{tabular}
\label{grgp3}
\end{equation}

With the help of these notations we can state $\left( \ref{r14}\right) $
simply as%
\begin{eqnarray}
&&\sum_{n=r}^{\infty }\binom{n}{r}H_{n}^{\left( \alpha \right) }\frac{r!}{%
n^{r}}f\left( n\right) x^{n}  \label{r14+} \\
&=&\frac{1}{\left( 1-z\right) ^{\alpha }}\sum_{n=r}^{\infty }\frac{f^{\left(
n\right) }\left( 0\right) }{n!}\left[ _{r}F_{n,\alpha }^{h}\left( \frac{x}{%
1-x}\right) -_{r}F_{n,\alpha }\left( \frac{x}{1-x}\right) \ln \left(
1-x\right) \right] .  \notag
\end{eqnarray}

\begin{remark}
Specializing $x=1$ in $\left( \ref{r15}\right) $\ we get "hyperharmonic $r-$
geometric numbers" as%
\begin{equation}
_{r}F_{n,\alpha }^{h}=\sum_{k=0}^{n}\QATOPD\{ \} {n}{k}_{r}k!H_{k}^{\left(
\alpha \right) }.  \label{r17}
\end{equation}
\end{remark}

The first few hyperharmonic $r-$geometric numbers are:

Case $\alpha =2$

\begin{equation}
\begin{tabular}{|l|l|l|}
\hline
$_{r}F_{n,2}^{h}$ & $r=1$ & $r=2$ \\ \hline
$n=0$ & $0$ & $0$ \\ \hline
$n=1$ & $1$ & $0$ \\ \hline
$n=2$ & $6$ & $5$ \\ \hline
$n=3$ & $42$ & $36$ \\ \hline
$n=4$ & $346$ & $304$ \\ \hline
\end{tabular}
\label{hhrfn2}
\end{equation}

Case $\alpha =3$

\begin{equation}
\begin{tabular}{|l|l|l|}
\hline
$_{r}F_{n,3}^{h}$ & $r=1$ & $r=2$ \\ \hline
$n=0$ & $0$ & $0$ \\ \hline
$n=1$ & $1$ & $0$ \\ \hline
$n=2$ & $8$ & $7$ \\ \hline
$n=3$ & $69$ & $61$ \\ \hline
$n=4$ & $674$ & $605$ \\ \hline
\end{tabular}
\label{hhrfn3}
\end{equation}

and specializing $x=1$ in $\left( \ref{r16}\right) $\ gives "general $r-$
geometric numbers" as%
\begin{equation}
_{r}F_{n,\alpha }=\frac{1}{\Gamma \left( \alpha \right) }\sum_{k=0}^{n}%
\QATOPD\{ \} {n}{k}_{r}\Gamma \left( k+\alpha \right) .  \label{r18}
\end{equation}

The first few general $r-$geometric numbers are:

Case $\alpha =2$

\begin{equation}
\begin{tabular}{|l|l|l|}
\hline
$_{r}F_{n,2}$ & $r=1$ & $r=2$ \\ \hline
$n=0$ & $0$ & $0$ \\ \hline
$n=1$ & $2$ & $0$ \\ \hline
$n=2$ & $8$ & $6$ \\ \hline
$n=3$ & $44$ & $36$ \\ \hline
$n=4$ & $308$ & $264$ \\ \hline
\end{tabular}
\label{grgn2}
\end{equation}%
and

Case $\alpha =3$

\begin{equation}
\begin{tabular}{|l|l|l|}
\hline
$_{r}F_{n,3}$ & $r=1$ & $r=2$ \\ \hline
$n=0$ & $0$ & $0$ \\ \hline
$n=1$ & $3$ & $0$ \\ \hline
$n=2$ & $15$ & $12$ \\ \hline
$n=3$ & $99$ & $84$ \\ \hline
$n=4$ & $807$ & $708$ \\ \hline
\end{tabular}
\label{grgn3}
\end{equation}

Thanks to the following corollary of Proposition $\left( \ref{phht}\right) $
we have closed forms of some series related to hyperharmonic numbers and
binomial coefficients.

\begin{corollary}
\begin{equation}
\sum_{n=r}^{\infty }\binom{n}{r}r!n^{m-r}H_{n}^{\left( \alpha \right) }x^{n}=%
\frac{1}{\left( 1-z\right) ^{\alpha }}\left[ _{r}F_{m,\alpha }^{h}\left( 
\frac{x}{1-x}\right) -_{r}F_{m,\alpha }\left( \frac{x}{1-x}\right) \ln
\left( 1-x\right) \right] .  \label{r19}
\end{equation}
\end{corollary}

\begin{proof}
For a positive integers $m\geq r,$ setting $f\left( x\right) =x^{m}$ in $%
\left( \ref{r14+}\right) $ gives $\left( \ref{r19}\right) $.
\end{proof}

\begin{remark}
Specializing the values of $r$, $m$ and $\alpha $ in $\left( \ref{r19}%
\right) $\ one can get closed forms of several hyperharmonic numbers series.
\end{remark}

Now we extend the formula $\left( \ref{gs5}\right) $\ to hyperharmonic
number series.

\begin{proposition}
\begin{eqnarray}
&&\sum_{n=r}^{\infty }\left( \sum_{k=0}^{n-r}\binom{k}{r}\binom{n+s-k}{s}%
r!k^{m-r}H_{k}^{\left( \alpha \right) }\right) x^{n}  \notag \\
&=&\sum_{n=r}^{\infty }\left( \sum_{0\leq k_{1}\leq k_{2}\leq \cdots \leq
k_{s+1}\leq n}\binom{k_{1}}{r}r!k_{1}^{m-r}H_{k_{1}}^{\left( \alpha \right)
}\right) x^{n}  \label{mhhrs} \\
&=&\frac{1}{\left( 1-x\right) ^{\alpha +s+1}}\left\{ _{r}F_{m}^{h}\left( 
\frac{x}{1-x}\right) -_{r}F_{m}\left( \frac{x}{1-x}\right) \ln \left(
1-x\right) \right\}  \notag
\end{eqnarray}
\end{proposition}

\begin{proof}
Multiplying both sides of the equation $\left( \ref{r19}\right) $ with the
Newton binomial series gives the statement.
\end{proof}

\begin{remark}
Also special values of $r$, $m$, $s$ and $\alpha $ in $\left( \ref{mhhrs}%
\right) $\ gives closed forms of several mutiplicative hyperharmonic numbers
series.
\end{remark}

\begin{remark}
Using $\left( \ref{24}\right) $ we get an alternative expression of
hyperharmonic $r-$ geometric polynomials and numbers as%
\begin{equation}
_{r}F_{n,\alpha }^{h}\left( x\right) =\sum_{k=0}^{n}\QATOPD\{ \} {n}{k}_{r}%
\QATOPD[ ] {k+\alpha }{\alpha +1}_{\alpha }x^{k},  \label{r20}
\end{equation}%
\begin{equation}
_{r}F_{n,\alpha }^{h}=\sum_{k=0}^{n}\QATOPD\{ \} {n}{k}_{r}\QATOPD[ ] {%
k+\alpha }{\alpha +1}_{\alpha }.  \label{r21}
\end{equation}
\end{remark}

\subsection{Hyperharmonic $r-$exponential polynomials and numbers}

\begin{definition}
For positive integers $m$ and $r$ "hyperharmonic $r-$exponential
polynomials" are defined as%
\begin{equation}
_{r}\phi _{n,\alpha }^{h}\left( x\right) =\sum_{k=0}^{n}\QATOPD\{ \}
{n}{k}_{r}H_{k}^{\left( \alpha \right) }x^{k}\text{.}  \label{r22}
\end{equation}
\end{definition}

The first few hyperharmonic $r-$exponential polynomials are:

Case $\alpha =2$

\begin{equation}
\begin{tabular}{|l|l|l|}
\hline
$_{r}\phi _{n,2}^{h}\left( x\right) $ & $r=1$ & $r=2$ \\ \hline
$n=0$ & $0$ & $0$ \\ \hline
$n=1$ & $x$ & $0$ \\ \hline
$n=2$ & $x+\frac{5}{2}x^{2}$ & $\frac{5}{2}x^{2}$ \\ \hline
$n=3$ & $x+\frac{15}{2}x^{2}+\frac{13}{3}x^{3}$ & $5x^{2}+\frac{13}{3}x^{3}$
\\ \hline
$n=4$ & $x+\frac{35}{2}x^{2}+26x^{3}+\frac{77}{12}x^{4}$ & $10x^{2}+\frac{65%
}{3}x^{3}+\frac{77}{12}x^{4}$ \\ \hline
\end{tabular}
\label{hhrep2}
\end{equation}

Case $\alpha =3$

\begin{equation}
\begin{tabular}{|l|l|l|}
\hline
$_{r}\phi _{n,3}^{h}\left( x\right) $ & $r=1$ & $r=2$ \\ \hline
$n=0$ & $0$ & $0$ \\ \hline
$n=1$ & $x$ & $0$ \\ \hline
$n=2$ & $x+\frac{7}{2}x^{2}$ & $\frac{7}{2}x^{2}$ \\ \hline
$n=3$ & $x+\frac{21}{2}x^{2}+\frac{47}{6}x^{3}$ & $7x^{2}+\frac{47}{6}x^{3}$
\\ \hline
$n=4$ & $x+\frac{49}{2}x^{2}+47x^{3}+\frac{171}{12}x^{4}$ & $14x^{2}+\frac{%
235}{6}x^{3}+\frac{171}{12}x^{4}$ \\ \hline
\end{tabular}
\label{hhrep3}
\end{equation}

Hence "hyperharmonic $r-$exponential numbers" are defined as%
\begin{equation}
_{r}\phi _{n,\alpha }^{h}=\sum_{k=0}^{n}\QATOPD\{ \} {n}{k}_{r}H_{k}^{\left(
\alpha \right) }.  \label{r23}
\end{equation}

The first few hyperharmonic $r-$exponential numbers are:

Case $\alpha =2$

\begin{equation}
\begin{tabular}{|l|l|l|}
\hline
$_{r}\phi _{n,2}^{h}$ & $r=1$ & $r=2$ \\ \hline
$n=0$ & $0$ & $0$ \\ \hline
$n=1$ & $1$ & $0$ \\ \hline
$n=2$ & $\frac{7}{2}$ & $\frac{5}{2}$ \\ \hline
$n=3$ & $\frac{77}{6}$ & $\frac{28}{3}$ \\ \hline
$n=4$ & $\frac{611}{12}$ & $\frac{457}{12}$ \\ \hline
\end{tabular}
\label{hhren2}
\end{equation}

Case $\alpha =3$

\begin{equation}
\begin{tabular}{|l|l|l|}
\hline
$_{r}\phi _{n,3}^{h}$ & $r=1$ & $r=2$ \\ \hline
$n=0$ & $0$ & $0$ \\ \hline
$n=1$ & $1$ & $0$ \\ \hline
$n=2$ & $\frac{9}{2}$ & $\frac{7}{2}$ \\ \hline
$n=3$ & $\frac{58}{3}$ & $\frac{89}{6}$ \\ \hline
$n=4$ & $\frac{347}{4}$ & $\frac{809}{12}$ \\ \hline
\end{tabular}
\label{hhren3}
\end{equation}

\begin{remark}
We extend the relation $\left( \ref{r11}\right) $ as
\end{remark}

\begin{equation*}
_{r}F_{n,\alpha }^{h}\left( z\right) =\int_{0r}^{\infty }\phi _{n,\alpha
}^{h}\left( z\lambda \right) e^{-\lambda }d\lambda .
\end{equation*}


\begin{thebibliography}{10}
\bibitem[1]{AS} Abramowitz, M and Stegun, I.\emph{\ Handbook of Mathematical
Functions with Formulas, Graphs, and Mathematical Tables}, 9th printing. New
York: Dover, p. 824, 1972.

\bibitem[2]{BL1} Bell, E. T. \emph{Exponential polynomials}, Annals of
Mathematics, vol. 35, no. 2, pp. 258--277, (1934).

\bibitem[3]{BL2} Bell, E. T. \emph{Exponential numbers}, Amer. Math. Monthly
41, 411-419, (1934).

\bibitem[4]{BG} Benjamin, A.T., Gaebler D. and Gaebler, R. \emph{A
combinatorial approach to hyperharmonic numbers}, Integers: Electron. J.
Combin. Number Theory 3 (2003), pp. 1--9 \#A15.

\bibitem[5]{BQ} Benjamin, A. T. and Quinn, J. J. \emph{Proofs that Really
Count: The Art of Combinatorial Proof}, MAA, 2003.

\bibitem[6]{B} Boyadzhiev, Khristo N. \emph{A Series transformation formula
and related polynomials}, In. J. Math. Math. Sc. 2005: 23 (2005), 3849-3866.

\bibitem[7]{B2} Boyadzhiev, Khristo N. \emph{Exponential Polynomials,
Stirling Numbers and Evaluation of some Gamma Integrals, }Abstract and
Applied Analysis, Volume 2009, Article ID 168672.

\bibitem[8]{Br} Broder, A. Z. \emph{The r-Stirling numbers}, Discrete Math.
49 (1984), 241-259.

\bibitem[9]{CA1} Carlitz, L. \emph{Weighted Stirling numbers of the first
and second kind-I}, The Fibonacci Quarterly, 18 (1980), 147-162.

\bibitem[10]{CA2} Carlitz, L. \emph{Weighted Stirling numbers of the first
and second kind-II}, The Fibonacci Quarterly, 18 (1980), 242-257.

\bibitem[11]{C} Comtet, L. \emph{Advanced Combinatorics. The Art of Finite
and Infinite Expansions}\textit{, Revised and enlarged edition}, D. Riedel
Publishing Co., Dordrecht, 1974.

\bibitem[12]{CG} Conway J. H. and Guy R. K., \emph{The book of numbers}, New
York, Springer-Verlag, 1996.

\bibitem[13]{DM} Dil, A. and Mezo, I. \emph{A Symmetric Algorithm for
Hyperharmonic and Fibonacci Numbers}, Applied Mathematics and Computation
206 (2008), 942--951.

\bibitem[14]{DK} Dil, A. and Kurt, V. \emph{Investigating Geometric and
exponential Polynomials with Euler-Seidel Algorithm, }submitted\emph{. }%
Available at http://arxiv.org/abs/0908.2585.

\bibitem[15]{DK1} Dil, A. and Kurt, V. \emph{Polynomials Related to Harmonic
Numbers and Evaluation of Harmonic Number Series I, }submitted\emph{.}

\bibitem[16]{G} Grunert, J. A. \emph{Uber die Summerung der Reihen...},
Journal f\"{u}r die reine und angewandte Mathematik, vol. 25, pp. 240--279,
1843.

\bibitem[17]{GKP} Graham, R. L., Knuth D. E. and Patashnik O., \emph{%
Concrete Mathematics}, Addison Wesley, 1993.

\bibitem[18]{M} Mezo, I.\emph{\ r-exponential numbers}, arXiv:0909.4417,
2009.

\bibitem[19]{MD} Mezo, I., Dil, A. \emph{Hyperharmonic series involving
Hurwitz zeta function, }Journal of Number Theory, 130, 2, 2010, 360-369.

\bibitem[20]{GN} Nyul, G. \emph{Ordered r-exponential numbers and r-Eulerian
numbers} (in Hungarian), Diophantine and Cryptography Days in Sopron, 11
October 2008, Sopron.

\bibitem[21]{Ri} Riordan, J. \emph{Combinatorial Analysis}, John Wiley, New
York, 1958.

\bibitem[22]{S} Schwatt, I. J. \emph{An Introduction to the Operations with
Series}, Chelsea, New York, 1962.

\bibitem[23]{ST} Tanny, Stephen M. \emph{On some numbers related to the
exponential numbers}, Canadian Mathematical Bulletin, Vol 17, (1974), No 5,
733-738.

\bibitem[24]{W} Wilf, Herbert. S. \emph{Generatingfunctionology, }Academic
Press, 1993.
\end{thebibliography}
\end{document}